\documentclass[12pt]{article}

\usepackage[leqno]{amsmath}
\usepackage{amsthm,amscd}
\usepackage{amsfonts} 
\usepackage{amssymb} 

\usepackage[dvips]{graphicx}

\newtheorem{theorem}{Theorem}[section]

\def\kmms{\kern-\mathsurround}

\newcommand{\calX}{{\cal X}}

\newcommand{\calZ}{{\cal Z}}

\newcommand{\bbR}{{\mathbb R}}
\newcommand{\bbS}{{\mathbb S}}

\begin{document}

\title{Estimation of the shape of 
the 
density contours 
of 
star-shaped distributions
}
\author{
Hidehiko Kamiya\footnote
{
{\rm This work was partially supported by 
JSPS 
KAKENHI Grant Number 25400201. } 
} 
\\ 
{\it Faculty of Economics, Osaka University of Economics}
}
\date{January 2017}

\maketitle

\begin{abstract} 
Elliptically contoured distributions generalize 
the multivariate normal distributions 
in such a way that 
the density generators need not be exponential. 
However, as the name suggests, 
elliptically contoured distributions 
remain to be restricted 
in that the similar 
density contours 
ought to be 
elliptical. 
Kamiya, Takemura and Kuriki 
[Star-shaped distributions and their generalizations, 
{\it Journal of Statistical Planning and Inference} {\bf 138} (2008), 3429--3447] 
proposed star-shaped distributions, 
for which 
the density contours 
are allowed to be 
boundaries of 
arbitrary 
similar 
star-shaped sets. 
In the present paper, we propose 
a nonparametric estimator of the shape of 
the 
density contours 
of star-shaped distributions, 
and prove its strong consistency with respect to 
the Hausdorff distance. 
We illustrate our estimator by simulation. 


\end{abstract}

\section{Introduction} 
\label{sec:introduction}

Elliptically contoured distributions generalize 
the multivariate normal distributions 
in such a way that 
the density generators need not be exponential 
(Fang and Zhang \cite{FZ90}). 
In this way, 
the class of elliptically contoured distributions 
includes, 
for example, 
distributions whose tails are heavier than those of the multivariate normal 
distributions. 
However, as the name suggests, 
elliptically contoured distributions 
remain to be restricted 
in that the similar 
density contours 
ought to be 
elliptical. 
Hence, in particular, 
no skewed distributions are  
members of this class. 
Skew-elliptical distributions (Genton \cite{Gen04}) allow 
skewness by introducing an extra parameter into 
elliptically contoured distributions. 

Kamiya, Takemura and Kuriki \cite{KTK08} proposed 
a flexible class of distributions called 
star-shaped distributions, 
for which 
the density contours 
are allowed to be 
boundaries of 
arbitrary 
similar 
star-shaped sets 
(see also \cite{TK96}, \cite{KT08}). 
Essentially the same idea can be found in  
$v$-spherical distributions 
by Fern\'{a}ndez, Osiewalski and Steel \cite{FOS95JASA} 
and 
center-similar distributions 
by Yang and Kotz \cite{YK03}. 
Skewness 
as well as heavy-tailedness 
is allowed in star-shaped distributions. 
Thus, 
besides (centrally, 
reflectively 
or in some other ways) symmetric distributions 
such as elliptically contoured distributions and 
$l_q$-spherical distributions, 
the class of star-shaped distributions 
also 
includes asymmetric distributions. 

Kamiya, Takemura and Kuriki \cite{KTK08} 
studied 
distributional properties of star-shaped distributions, 
including 
independence of 
the ``length'' and the ``direction,'' 
and robustness of the distribution of the ``direction.''  
However, they did not 
explore 
inferential 
aspects of 
star-shaped distributions. 
From the perspective 
of 
\cite{KTK08}, 
the most important problem 
in 
the inference for 
star-shaped distributions is 
the estimation of the shape of 
the 
density contours.  

In the present paper, we propose 
a nonparametric estimator of the shape of 
the density contours. 
The point is that the density of the usual 
direction 
under a star-shaped distribution is in one-to-one 
correspondence with a function which determines the shape 
of the density contours. 
Thus, by nonparametrically estimating 
the density of the direction, 
we can obtain a nonparametric estimator of the shape. 
We prove its strong consistency with respect to 
the Hausdorff distance. 

In a recent paper, Liebscher and Richter \cite{LR16} 
presented 
examples of 
parametric modeling and estimation 
concerning 
the shape of the density contours 
of two-dimensional star-shaped distributions 
(Section 2.2 as well as Sections 3.3 and 3.4 of \cite{LR16}). 
They 
also investigated 
estimation about 
many 
other aspects of 
star-shaped distributions.  

The organization of this paper is as follows. 
We describe a star-shaped distribution and define the shape of its 
density contours in Section \ref{sec:star-shaped}. 
Next, we propose an estimator of 
the shape of 
the 
density contours of a star-shaped distribution 
in Section \ref{sec:estimator}, 
and prove 
its strong consistency in Section \ref{sec:consistency}. 
We illustrate our estimator by simulation 
in Section \ref{sec:simulation}, 
and conclude with some remarks 
in Section \ref{sec:concluding-remarks}. 

\section{Star-shaped distribution and the 
shape of its density contours}  
\label{sec:star-shaped} 

In this section, we describe a star-shaped distribution 
and define the shape of its density contours. 

Suppose a random vector 
$x \in \calX := \bbR^p \setminus \{ {\bf 0} \}, 
\ p \ge 2$, 
is distributed as 
\begin{equation} 
\label{eq:h(r(x))} 
x \sim h(r(x))dx, 
\end{equation} 
where $r: 
\calX 
\to \bbR_{>0}$ is continuous 
and 
equivariant under the action of 
the positive real numbers: 
$r(cx)=cr(x)$ for all $c \in \bbR_{>0}$. 
In \eqref{eq:h(r(x))}, it is implicitly assumed that 
$h: \bbR_{>0}\to \bbR_{\ge 0}$ satisfies 
$0<\int_0^{\infty}h(r)r^{p-1}dr<\infty$. 
In the particular case that $r(x)=(x^T\Sigma^{-1} x)^{1/2}$ 
for a positive definite matrix $\Sigma$ 
($x^T$ denotes 
the 
transpose of the column vector $x$) 
and  
that the density generator $h((-2( \, \cdot \, ))^{1/2})$ is 
exponential: 
$h((-2( \, \cdot \, ))^{1/2})\propto \exp( \, \cdot \, )$, 
we obtain the multivariate normal distribution $N_p({\bf 0}, \Sigma)$. 

Define 
\begin{equation} 
\label{eq:Z} 
Z:=\left\{ x \in \calX: 
r(x)=1 \right\}, 
\end{equation}  
and write $cZ:=\{ cz: z \in Z \}$ for $c \in \bbR_{\ge 0}$. 
Then the density $h(r(x))$ is constant on each of $cZ \subset \calX, \ c \in \bbR_{>0}$: $h(r(x))=h(c)$ for all $x \in cZ$. 
When 
$h: \bbR_{>0} \to \bbR_{\ge 0}$ is injective 
(e.g., strictly 
decreasing), 
each $cZ, \ c \in \bbR_{>0}$, 
is 
a contour of the density $h(r(x))$: 
$cZ=\{ x \in \calX: h(r(x))=h(c)\}$, 
but in general, a contour of the density is 
a 
union of some $cZ$'s: 
$\{ x \in \calX: h(r(x))=t \}=\bigcup_{c\in h^{-1}(\{ t \})}cZ, 
\ t \in \bbR_{\ge 0}$.  

Noticing that 
$\calZ:=\bigcup_{0 \le c \le 1}cZ 
\subset \calX \cup \{ {\bf 0} \}
=\bbR^p$ 
is a star-shaped set with respect to the origin, 
we say that $x$ in \eqref{eq:h(r(x))} has a star-shaped distribution. 
Also, we call $Z$ the shape 
of the density contours of 
this 
star-shaped distribution, 
including cases where 
$h$ is not injective. 
When $h$ is strictly decreasing, $\calZ$ is 
a density level set: 
$\calZ=\{ x \in \calX : h(r(x)) \ge h(1)\}\cup \{ {\bf 0} \}$. 

The focus of this paper is the estimation of 
the shape 
$Z$ in \eqref{eq:Z}. 
In the next section, 
we propose an estimator of the form 
\[ 
\left\{ \hat{f}_n(u)^{\frac{1}{p}} u: u \in \bbS^{p-1} \right\}, 
\] 
where $\bbS^{p-1}$ is the unit sphere in $\bbR^p$ 
and $\hat{f}_n(u)$ is a directional 
density estimator 
based on the directions of a sample from  \eqref{eq:h(r(x))}.  

\section{Estimation of the shape}  
\label{sec:estimation} 

In this section, 
we propose an estimator of 
the shape of 
the 
density contours of a star-shaped distribution 
(Section \ref{sec:estimator}), 
and prove 
its strong consistency (Section \ref{sec:consistency}). 

\subsection{Proposed estimator} 
\label{sec:estimator} 

In this subsection, we propose an estimator of the shape $Z$. 

Let $\| \cdot \|$ denote the Euclidean norm. 
Under \eqref{eq:h(r(x))}, the direction $u:=x/\| x \| 
\in \bbS^{p-1}$ 
is distributed as 
\begin{equation} 
\label{eq:u-fdu-f=}  
u \sim f(u)du
\text{ \quad with \ } f(u):=c_0r(u)^{-p}, 
\end{equation} 
where $du$ stands for the volume 
element of $\bbS^{p-1}$ and 
$c_0=1/\int_{\bbS^{p-1}}r(u)^{-p}du
=\int_0^{\infty}h(r)r^{p-1}dr$ 
(Theorem 4.1 of \cite{KTK08}). 
Note 
the function $f: \bbS^{p-1} \to \bbR_{\ge 0}$ 
in \eqref{eq:u-fdu-f=} is continuous 
and satisfies 
$f(u)>0$ for all $u \in \bbS^{p-1}$.  
Throughout this section (Section \ref{sec:estimation}), 
we assume $r( \, \cdot \, )$ is taken so that 
$\int_{\bbS^{p-1}}r(u)^{-p}du=1$ and hence $c_0=1$. 

Now, we can write $r(u)=f(u)^{-1/p}$ for $u \in \bbS^{p-1}$.  
Thus, for 
$x \in \calX$, 
the condition that $r(x)=1$ is equivalent to $\| x \| = 1/r(x/\| x \|)=
f(x/\| x \|)^{1/p}$. 
Hence 
$Z=\{ x \in \calX: 
r(x)=1 \}
=\{ f(u)^{1/p} u: u \in \bbS^{p-1} \}$,  
and we can estimate $Z$ by estimating 
the density $f(u)$ of $u=x/\| x \|$.  

Suppose we are given 
an i.i.d. sample 
$x_1,\ldots,x_n$ from \eqref{eq:h(r(x))}, 
and consider estimating 
$f(u)$ based on $u_1,\ldots,u_n$, 
where $u_i:=x_i/\| x_i \|, \ i=1,\ldots,n$. 

Let $\hat{f}_n(u)$ be 
an estimator of $f(u)$ such that $\hat{f}_n(u)\ge 0$ for 
all $u \in \bbS^{p-1}$. 
Define the estimator $\hat{Z}_n$ of $Z$ 
by 
\[
\hat{Z}_n:=\left\{ \hat{f}_n(u)^{\frac{1}{p}} u: u \in \bbS^{p-1} \right\}. 
\] 
Then $\hat{\calZ}_n:=\bigcup_{0 \le c \le 1}c\hat{Z}_n$ is also a star-shaped set 
with respect to the origin. 

\subsection{Strong consistency} 
\label{sec:consistency} 

In this subsection, we 
prove 
strong consistency of our estimator $\hat{Z}_n$.  

Let $\delta_H(\hat{Z}_n, Z)$ be the Hausdorff distance 
between $\hat{Z}_n$ and $Z$: 
\[ 
\delta_H(\hat{Z}_n, Z)
:= \inf 
\left\{ \delta>0: 
\hat{Z}_n \subset Z + B(\delta), \  
Z \subset \hat{Z}_n + B(\delta)
\right\}, 
\] 
where $B(\delta):=\{ x \in \bbR^p: \| x \| \le \delta\}$, and 
$+$ denotes the Minkowski sum. 
Similarly, let 
$\delta_H(\hat{\calZ}_n, \calZ)
= \inf 
\{ \delta>0: 
\hat{\calZ}_n \subset \calZ + B(\delta), \  
\calZ \subset \hat{\calZ}_n + B(\delta)
\}$ 
be the Hausdorff distance between 
$\hat{\calZ}_n$ and $\calZ$. 
We note that $\hat{Z}_n$ and $\hat{\calZ}_n$ 
may not be compact.  
The purpose of this section is to show that, 
under some conditions, 
$\delta_H(\hat{Z}_n, Z)$ 
and 
$\delta_H(\hat{\calZ}_n, \calZ)$ 
converge to 
zero 
almost surely. 

We begin by proving that  
$\delta_H(\hat{Z}_n, Z)$ 
and 
$\delta_H(\hat{\calZ}_n, \calZ)$
are 
bounded by 
$d_n:=
\sup_{u \in \bbS^{p-1}} \vert \hat{f}_n(u)^{1/p}-f(u)^{1/p} \vert$: 
\begin{equation} 
\label{eq:d_H<=} 
\delta_H(\hat{Z}_n, Z) 
\le 
d_n, 
\qquad 
\delta_H(\hat{\calZ}_n, \calZ) 
\le 
d_n. 
\end{equation}  
Let 
$z_0
=\tilde{c}_0 f(u_0)^{1/p}u_0$ $(0 \le \tilde{c}_0 \le 1, \ u_0 \in \bbS^{p-1})$ be 
an arbitrary point of $\calZ$. 
Take $z'_0=\tilde{c}_0\hat{f}_n(u_0)^{1/p}u_0 \in \hat{\calZ}_n$. 
Then $\| z'_0-z_0 \| 
= \tilde{c}_0 \vert \hat{f}_n(u_0)^{1/p}-f(u_0)^{1/p} \vert 
\le d_n$, 
and thus 
$z_0 \in 
\hat{\calZ}_n + B(d_n)$. 
This argument implies that $\calZ \subset \hat{\calZ}_n + B(d_n)$. 
Similarly, $\hat{\calZ}_n \subset \calZ + B(d_n)$ holds true. 
Therefore, the second inequality 
in 
\eqref{eq:d_H<=} is proved. 
The proof of the first inequality 
in 
\eqref{eq:d_H<=} is similar. 

Next we want to verify that 
$
d_n 
\to 0$ 
almost surely 
for estimators $\hat{f}_n(u)$ 
having a certain property. 

For each $u\in \bbS^{p-1}$ and each $n$, we can write 
\begin{equation} 
\label{eq:linear} 
\hat{f}_n(u)^{\frac{1}{p}}
=f(u)^{\frac{1}{p}}+\frac{1}{p}
f^*_n(u)^{\frac{1}{p}-1} \left( \hat{f}_n(u)-f(u)\right) 
\end{equation} 
for some $f^*_n(u)$ between 
$\hat{f}_n(u)$ and $f(u)$. 

Let $\epsilon_n:=\sup_{u \in \bbS^{p-1}}\vert \hat{f}_n(u)-f(u)\vert$. 
Then we have $f^*_n(u) \ge f(u)-\epsilon_n$ for 
all $u \in \bbS^{p-1}$ and all $n$, 
and thus 
\begin{equation} 
\label{eq:inf-f*>=} 
\inf_{u \in \bbS^{p-1}} f^*_n(u) 
\ge 
\inf_{u \in \bbS^{p-1}} f(u)
- \epsilon_n
\end{equation} 
for all $n$. 
Since $f: \bbS^{p-1} \to \bbR_{\ge 0}$ 
is continuous, 
$\bbS^{p-1}$ is compact  
and $f(u)>0$ for all $u\in \bbS^{p-1}$,  
we have $c_f:=\inf_{u \in \bbS^{p-1}}f(u)
=\min_{u \in \bbS^{p-1}}f(u)>0$. 
Now, suppose the estimator $\hat{f}_n(u)$ satisfies 
\begin{equation} 
\label{eq:en-goes-to-0} 
\epsilon_n
=\sup_{u \in \bbS^{p-1}}
\left\vert \hat{f}_n(u)-f(u) \right\vert 
\to 0 
\quad {\rm a.s. } 
\end{equation}  
Then, with probability one, 
we have 
$\epsilon_n<c_f/2$ for all sufficiently large $n$. 
Together with this fact, 
inequality 
\eqref{eq:inf-f*>=} implies that, 
with probability one, 
\begin{equation} 
\label{eq:inf-f*>0} 
\inf_{u \in \bbS^{p-1}} f^*_n(u) 
\ge 
c_f-\epsilon_n 
> \frac{c_f}{2} 
\end{equation} 
for all sufficiently large $n$. 

It follows from 
\eqref{eq:linear} and \eqref{eq:inf-f*>0} 
that, 
with probability one, 
\begin{eqnarray*} 
d_n
=\sup_{u \in \bbS^{p-1}}
\left\vert
\hat{f}_n(u)^{\frac{1}{p}}-f(u)^{\frac{1}{p}} 
\right\vert 
&\le& \frac{1}{p} \left\{ \inf_{u \in \bbS^{p-1}}f^*_n(u)\right\}^{\frac{1}{p}-1}  
\cdot \sup_{u \in \bbS^{p-1}}\left\vert \hat{f}_n(u)-f(u) \right\vert \\ 
&\le& \frac{1}{p} \left( \frac{c_f}{2} \right)^{\frac{1}{p}-1} 
\epsilon_n 
\end{eqnarray*} 
for all sufficiently large $n$.  
Therefore, by \eqref{eq:en-goes-to-0} 
we obtain 
$
d_n 
\to 0 \ {\rm a.s.}$, 
as was to be verified.  

Now, 
for estimating a general density $f(u)$ on 
$\bbS^{p-1}, \ p \ge 2$ 
(i.e., not necessarily $f(u)$ in \eqref{eq:u-fdu-f=}) 
based on an i.i.d. sample $u_1,\ldots, u_n$ 
from $f(u)du$,  
we can use the following kernel density  
estimator  
(Hall, Watson and Cabrera \cite{HWC87}, Bai, Rao and Zhao \cite{BRZ88}): 
\begin{equation} 
\label{eq:def-of-hat-f} 
\hat{f}_n(u)
=\frac{C(\eta)}{n\eta^{p-1}}
\sum_{i=1}^n L \left( \frac{1-u^T u_i}{\eta^2} \right), 
\quad u \in \bbS^{p-1}, 
\end{equation} 
where 
$\eta=\eta_n>0$, 
$C(\eta):=\eta^{p-1}/\int_{\bbS^{p-1}}L((1-u^Ty)/\eta^2)du>0$ 
($y\in \bbS^{p-1}$), 
and 
$L: \bbR_{\ge 0}\to \bbR_{\ge 0}$ satisfies 
$0<\int_{0}^{\infty} L(s)s^{(p-3)/2}ds<\infty$.  
Notice that $C(\eta)$ does not depend on 
$y$ 
and can be written as  
$C(\eta)
=\eta^{p-1}/
\{ 
\omega_{p-1} 
\int_{-1}^{1}L((1-t)/\eta^2)(1-t^2)^{(p-3)/2}dt 
\}
=1/
\{ 
\omega_{p-1}
\int_0^{2/\eta^2}L(s)s^{(p-3)/2}(2-\eta^2 s)^{(p-3)/2}ds
\}, \ \omega_{p-1}:=2\pi^{(p-1)/2}/\Gamma((p-1)/2)$ 
(equation (2$\cdot$2) 
of 
\cite{HWC87}, 
equation (1.6) 
of 
\cite{BRZ88}). 
Recall, in passing, that 
the class of kernel estimators of the form 
\eqref{eq:def-of-hat-f} virtually ``contains asymptotically''  
the class of kernel estimators of the form 
$\tilde{f}_n(u)=(c_0(\kappa)/n)\sum_{i=1}^nK(\kappa u^Tu_i)$, 
$c_0(\kappa)=1/\int_{\bbS^{p-1}}K(\kappa u^Ty)du 
\ (y\in \bbS^{p-1})$, 
for a kernel $K$ and a smoothing parameter $\kappa>0$ 
(see Hall, Watson and Cabrera \cite{HWC87}). 
The choice $L(s)=\exp(-s)$, $K(s)=\exp(s)$ 
is the von Mises kernel. 

A sufficient condition for 
$\sup_{u \in \bbS^{p-1}}\vert \hat{f}_n(u)-f(u) \vert  
\to 0 \ {\rm a.s. }$ 
for a general density $f(u)$ on $\bbS^{p-1}$, 
$p \ge 2$, 
and its kernel estimator 
$\hat{f}_n(u)$ in \eqref{eq:def-of-hat-f} 
was obtained by Bai, Rao and Zhao \cite{BRZ88}, Theorem 2: 
$\sup_{u \in \bbS^{p-1}}\vert \hat{f}_n(u)-f(u) \vert  
\to 0 
\ {\rm a.s. } 
$  
holds true if the following 
conditions are satisfied: 
1. $f: \bbS^{p-1} \to \bbR_{\ge 0}$ is continuous; 
2. $L: \bbR_{\ge 0} \to \bbR_{\ge 0}$ is bounded;  
3. $L: \bbR_{\ge 0} \to \bbR_{\ge 0}$ is 
Riemann integrable on any finite interval 
in $\bbR_{\ge 0}$ 
with 
$\int_{0}^{\infty}
\sup_{t: \, \vert \sqrt{t}-\sqrt{s}\vert <1}L(t) 
\cdot 
s^{(p-3)/2}ds < \infty$; 
4. $\eta_n \to 0$ as $n\to \infty$;  
5. $n \eta_n^{p-1}/\log n\to \infty$ as $n\to \infty$. 

Note that 
under 
the fourth condition $\eta_n \to 0 \ (n \to \infty)$, 
we have 
$\lim_{n\to \infty}C(\eta_n)
=1/
\{ 
2^{(p-3)/2}\omega_{p-1}
\int_0^{\infty}L(s)s^{(p-3)/2}ds
\}
$ 
(equation (1.7) of \cite{BRZ88}).  

The preceding arguments yield the following result:  

\begin{theorem} 
\label{th:main} 
Let $x_1,\ldots , x_n \in \calX=\bbR^p\setminus \{ {\bf 0}\}, 
\ p \ge 2$, 
be an i.i.d. sample from 
a star-shaped distribution $h(r(x))dx$.  
Let 
$\hat{f}_n(u)
=(C(\eta)/(n\eta^{p-1}))
\sum_{i=1}^n L((1-u^T u_i)/\eta^2)$ be 
a kernel estimator of the density $f(u)$ of 
$u=x/\| x\|\in \bbS^{p-1}, \ x \sim h(r(x))dx$, 
based on $u_i=x_i/\| x_i\|, \ i=1,\ldots , n$. 

Assume 
the equivariant 
function 
$r: \calX \to \bbR_{>0}$ 
under the action of the positive real numbers 
is continuous 
and normalized so that 
$\int_{\bbS^{p-1}}r(u)^{-p}du=1$, 
and that 
$L: \bbR_{\ge 0} \to \bbR_{\ge 0}$ is bounded 
and satisfies 
$\int_{0}^{\infty}
\sup_{t: \, \vert \sqrt{t}-\sqrt{s}\vert <1}L(t) 
\cdot 
s^{(p-3)/2}ds < \infty$. 
Moreover, suppose 
$\eta=\eta_n>0$ is 
taken in such a way that 
$\eta_n \to 0$ and $n \eta_n^{p-1}/\log n\to \infty$ as $n\to \infty$.  

Then, 
$\hat{Z}_n=\{ \hat{f}_n(u)^{1/p}u: u \in \bbS^{p-1}\}$ 
is a strongly consistent estimator of 
the shape $Z=\{ x \in \calX: r(x)=1\}$ 
of the 
density contours of 
the 
star-shaped distribution $h(r(x))dx$ 
in the sense that 
the Hausdorff distance $\delta_H(\hat{Z}_n, Z)$  
between $\hat{Z}_n$ and $Z$ satisfies 
\[ 
\delta_H(\hat{Z}_n, Z) \to 0 
\quad {\rm a.s. } 
\] 
In addition, 
$\hat{\calZ}_n=\bigcup_{0 \le c \le 1}c\hat{Z}_n$ 
is a strongly consistent estimator of 
$\calZ=\bigcup_{0 \le c \le 1}cZ$: 
$\delta_H(\hat{\calZ}_n, \calZ) \to 0 \ {\rm a.s. }$ 

\end{theorem} 

It can easily be seen that $L(s)=e^{-s}$ 
and $L(s)=1(s<1)$ 
($= 1$ if $s<1$ and $0$ otherwise) 
satisfy 
$\int_{0}^{\infty}
\sup_{t: \, \vert \sqrt{t}-\sqrt{s}\vert <1}L(t) 
\cdot 
s^{(p-3)/2}ds < \infty$ 
and the other conditions 
of Theorem \ref{th:main}. 

\section{Illustrations by simulation} 
\label{sec:simulation} 

In this section, 
we illustrate our estimator by simulation.   

We consider star-shaped distributions in $\bbR^2$ 
and treat two shapes; 
one is the triangle in Examples 1.1 and 3.1 of 
Takemura and Kuriki 
\cite{TK96} (Section \ref{sec:triangle}), 
and the other is the unit $l_{1/2}$-sphere (Section \ref{sec:l_{1/2}}). 

In both cases, we use the von Mises kernel $L(s)=\exp(-s)$. 
We do not normalize $r( \, \cdot \, )$, 
so $c_0$ is not equal to one in general and 
our estimator of 
$Z=\{ x \in \calX: r(x)=1 \}$ is $\hat{Z}_n=\{ (\hat{f}_n(u)/c_0)^{1/2}u: u \in \bbS^1 \}$. 
We obtain the kernel estimator $\hat{f}_n(u)$ by making use of 
the 
R package \texttt{circular}\footnote{C. Agostinelli and U. Lund (2013). 
{R} package \texttt{circular}: Circular Statistics (version 0.4-7). 
URL  https://r-forge.r-project.org/projects/circular/ }. 
We select the bandwidth $1/\eta^2$ 
by simple trial and error. 
(If we did not know the true shape, we could use, e.g., cross-validation for 
minimizing the squared-error loss or the Kullback-Leibler loss 
in order to select the bandwidths 
(\cite{HWC87}).) 

Although we employ specific functions for $h( \, \cdot \, )$ below, 
these choices do not affect the estimation of $f(u)$ 
(and hence of $Z$)  
based on $u_1,\ldots,u_n$. 
This is because $u_i=x_i/\| x_i \|=z_i/\| z_i \|$ for 
$z_i:=x_i/r(x_i)\in Z, \ i=1, \ldots, n$, 
and the distribution of $z_i$ does not depend on $h( \, \cdot \, )$ 
(Theorem 4.1 of \cite{KTK08}). 

\subsection{Triangular shape} 
\label{sec:triangle} 

As in Examples 1.1 and 3.1 of \cite{TK96}, 
we take $r(x)=\max \{ -x(1), -x(2), x(1)+x(2) \}$ for 
$x=(x(1), x(2))$.  
Then the shape $Z$ is the triangle with vertices $P(2, -1)$, $Q(-1, 2)$ 
and $R(-1, -1)$.  
As is calculated in Example 3.1 of \cite{TK96}, 
we have 
$c_0
=1/\int_0^{2\pi} r(\cos \theta, \sin \theta)^{-2}d\theta 
=1/9$. 

Essentially as in Example 3.1 of \cite{TK96}, we choose $h(r)\propto \exp(-r^2/2)$, 
which necessarily implies $h(r)=c_0\exp(-r^2/2)=(1/9)\exp(-r^2/2)$ 
because of $c_0=\int_0^{\infty}h(r)rdr$ 
and $\int_0^{\infty}\exp(-r^2/2)rdr=1$. 
Hence, our star-shaped distribution is $(1/9)\exp\{-r(x)^2/2\}dx$. 

We can generate $x\sim (1/9)\exp\{-r(x)^2/2\}dx$ by $x=rz$, 
where $r \in \bbR_{>0}$ is distributed as the Rayleigh distribution 
with scale parameter 1 (i.e., $r^2\sim \chi^2(2)$), 
$z \in Z$ has density (with respect to the line element) 
$1/(9 \sqrt{2}), 1/9, 1/9$ on sides $PQ, QR, RP$, respectively 
(Example 3.1 of \cite{TK96}), 
and $r$ and $z$ are independently distributed.  

Our estimator of 
$Z$ is $\hat{Z}_n
=\{ (\hat{f}_n(u)/(1/9))^{1/2}u: u \in \bbS^1 \} 
=\{ 3\hat{f}_n(u)^{1/2}u: u \in \bbS^1 \}$. 

\begin{figure}[htbp]
 \begin{center}
\includegraphics*[width=.95\textwidth]{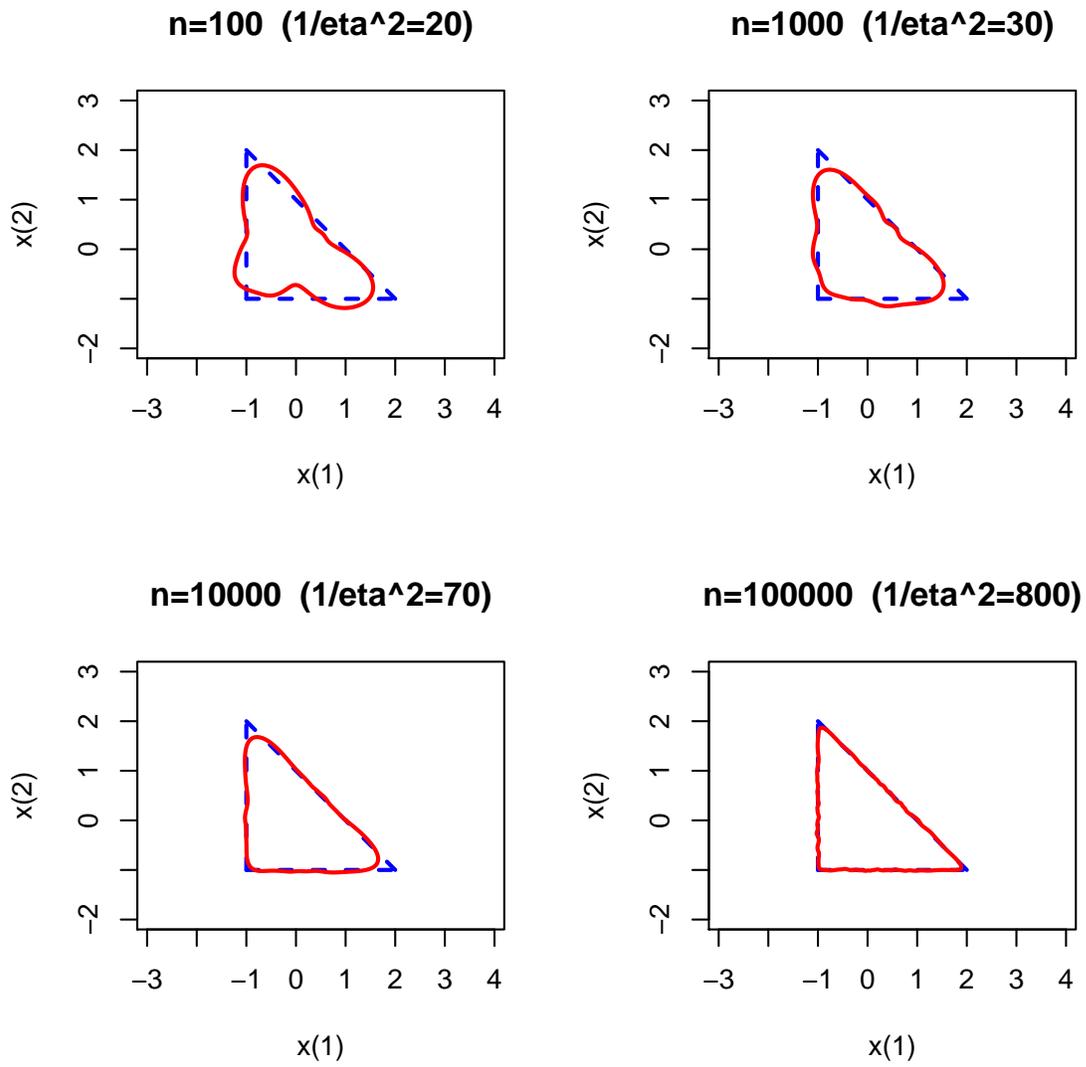}
\caption{Estimation of triangular shape.}
\label{fig:triangle}
 \end{center}
\end{figure}

The true shape $Z$ (blue, dashed line) 
and its estimator $\hat{Z}_n$ (red, solid line) 
for $n=100, 1000, 10000, 100000$ are shown in Figure \ref{fig:triangle}. 

\subsection{$l_{1/2}$-spherical shape} 
\label{sec:l_{1/2}} 

We take $r(x)=(\vert x(1) \vert^{1/2}+\vert x(2)\vert^{1/2})^2$, 
$x=(x(1), x(2))$.  
Then $Z=\{ (x(1), x(2)): \vert x(1) \vert^{1/2}+\vert x(2)\vert^{1/2}=1\}$ is 
the unit $l_{1/2}$-sphere. 
We can calculate  
$1/c_0 
=\int_{\bbS^1}r(u)^{-2}du 
=4\int_0^{\pi/2} \{ (\cos \theta)^{1/2}+ (\sin \theta)^{1/2} \}^{-4}d\theta 
=4/3$. 

We choose $h(r)\propto \exp(-2r^{1/2})$, so 
$h(r)=c_1\exp(-2r^{1/2})$, say. 
Then $c_0=\int_0^{\infty}h(r)rdr=c_1\int_0^{\infty}\exp(-2r^{1/2})rdr=(3/4)c_1$ 
and hence $c_1=(4/3)c_0=1$. 
Thus, our star-shaped distribution is $\exp\{-2r(x)^{1/2}\}dx$. 

This star-shaped distribution $\exp\{-2r(x)^{1/2}\}dx$ 
is obtained as 
the distribution of $x=(x(1), x(2))$ with $x(1)$ and $x(2)$ being  
independently distributed according to the $p$-generalized normal distribution 
with $p=1/2$ (this $p$ does not indicate the dimension of $\calX=\bbR^p\setminus \{ {\bf 0} \}$): 
$x(j) \sim \exp(-2\vert x(j) \vert^{1/2}), \ j=1, 2$. 
We generate $x(j), \ j=1,2$, by using 
the 
R package \texttt{pgnorm}\footnote{Steve Kalke (2015). 
\texttt{pgnorm}: The $p$-Generalized Normal Distribution. 
R package version 2.0. 
https://CRAN.R-project.org/package=pgnorm
}. 

Our estimator of 
$Z$ is 
$\hat{Z}_n
=\{ (2/\sqrt{3})\hat{f}_n(u)^{1/2}u: u \in \bbS^1 \}$. 

\begin{figure}[htbp]
 \begin{center}
\includegraphics*[width=.95\textwidth]{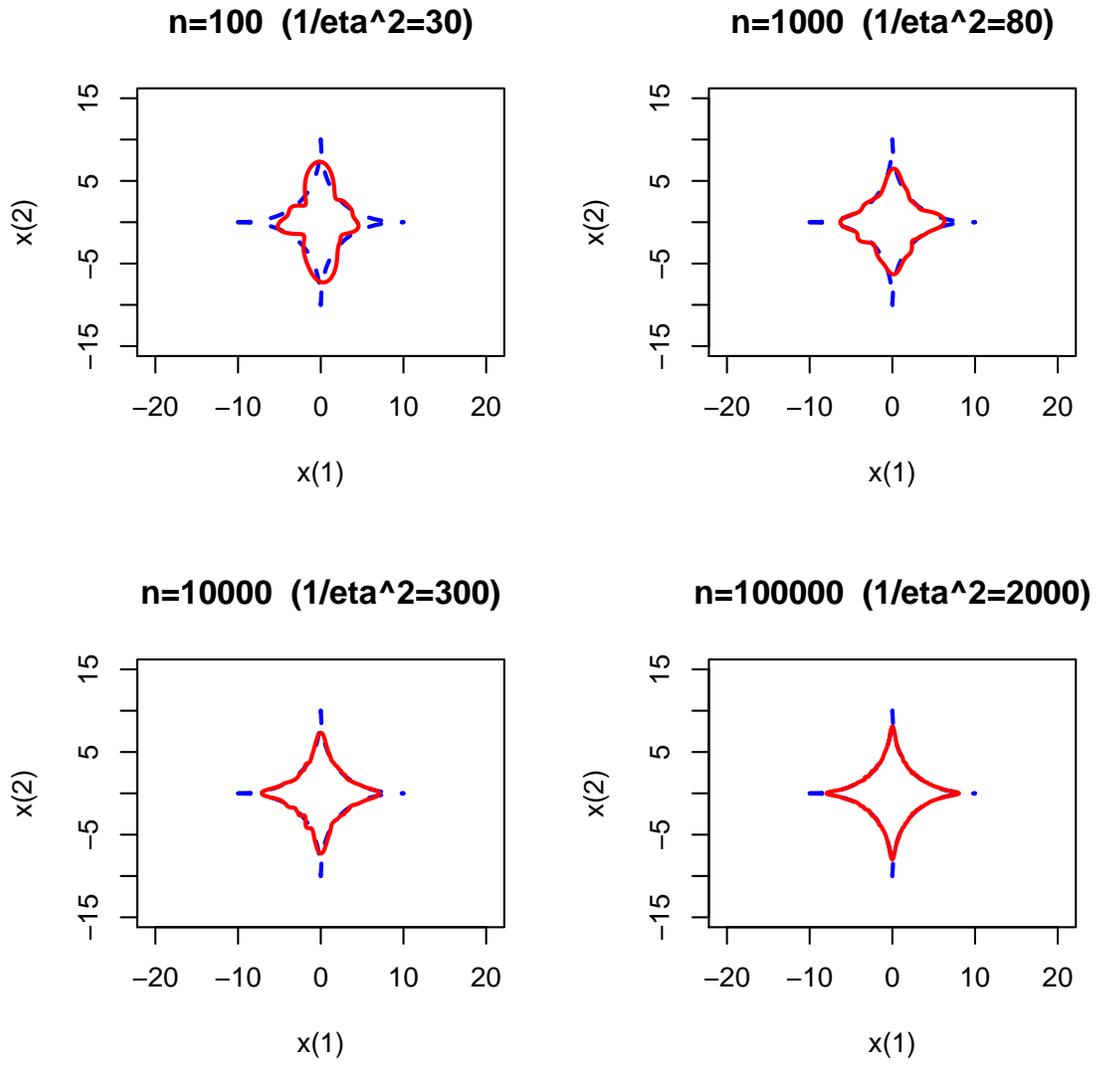}
\caption{Estimation of $l_{1/2}$-spherical shape.}
\label{fig:l-sphere}
 \end{center}
\end{figure}

For visibility, 
we enlarge the shape and its estimator, and display 
$10Z$ (blue, dashed line) 
and $10\hat{Z}_n$ (red, solid line) 
for $n=100, 1000, 10000, 100000$ in Figure \ref{fig:l-sphere}. 

\section{Concluding remarks}
\label{sec:concluding-remarks} 

In this paper, we proposed  
a nonparametric estimator of the shape of 
the 
density contours 
of star-shaped distributions, 
and proved its strong consistency with respect to 
the Hausdorff distance. 

We can introduce 
the 
location parameter and consider a star-shaped 
distribution whose density contours are 
(unions of) 
boundaries of 
star-shaped 
sets 
with respect to the location. 
In that case, one possibility for estimating the shape is 
to plug in an estimator of the location 
and use our proposed nonparametric estimator of the shape. 
We might be able to estimate the location by characterizing it in some way. 
For example, if the star-shaped distribution may be assumed to be 
centrally symmetric about 
the location and 
have a finite first moment, 
the location can be characterized as the mean and may be estimated by, 
e.g., the sample mean. 
If, instead, 
$h$ 
in \eqref{eq:h(r(x))} is 
strictly 
decreasing, 
the location can be regarded as the mode and 
be estimated by means of various methods 
for estimating the multivariate mode.

\end{document}